\documentclass[cleveref,a4paper,UKenglish]{lipics-v2021}
\usepackage{cite} 
\usepackage{amsmath}
\usepackage{amssymb}  
\usepackage{xspace}

\usepackage{tikz}
\usepackage{todonotes}
\tikzstyle{every picture} = [>=latex]
\usetikzlibrary{calc,arrows,decorations.markings,quotes}
\usetikzlibrary{backgrounds}

\def\ca#1{{\cal#1}}

\hideLIPIcs
\nolinenumbers

\title{Twin-width of Planar Graphs; a Short Proof}

\author{Petr Hlin{\v e}n\'y}{Masaryk University, Brno,
  Czech republic}{hlineny@fi.muni.cz}{https://orcid.org/0000-0003-2125-1514}{}
  
\authorrunning{P.\ Hlin\v{e}n\'y}
\Copyright{Petr Hlin\v{e}n\'y}



\begin{document}
\maketitle

\begin{abstract}
The fascinating question of the maximum value of twin-width on planar graphs
is nowadays not far from the final resolution; 
there is a lower bound of~$7$ coming from a construction by Kr\'al' and Lamaison [arXiv, September 2022],
and an upper bound of~$8$ by Hlin\v{e}n\'y and Jedelsk\'y [arXiv, October 2022].
The upper bound (currently best) of~$8$, however, is rather complicated and involved.
In the paper we give a short and simple self-contained proof that the twin-width of planar graphs is at most~$11$.
We believe that this short proof can also shed more light on the topic of upper bound(s)
on the twin-width of planar and beyond-planar graphs in general.
\end{abstract}

\section{Introduction}

The structural parameter twin-width was introduced in 2020 by Bonnet, Kim, Thomass\'e and Watrigant~\cite{DBLP:journals/jacm/BonnetKTW22}.
We consider it only for simple graphs (instead of general binary relational structures).

A \emph{trigraph} is a simple graph $G$ in which some edges are marked as {\em red}, and with respect to the red edges only, 
we naturally speak about \emph{red neighbours} and \emph{red degree} in~$G$. 
However, when speaking about edges, neighbours and/or subgraphs without further specification, we count both ordinary and red edges together as one edge set.
The edges of $G$ which are not red are sometimes called (and depicted) black for distinction.
For a pair of (possibly not adjacent) vertices $x_1,x_2\in V(G)$, we define a \emph{contraction} of the pair $x_1,x_2$ as the operation
creating a trigraph $G'$ which is the same as $G$ except that $x_1,x_2$ are replaced with a new vertex $x_0$ (said to {\em stem from $x_1,x_2$}) such that:
\begin{itemize}
\item 
the (full) neighbourhood of $x_0$ in $G'$ (i.e., including the red neighbours), denoted by $N_{G'}(x_0)$,
equals the union of the neighbourhoods $N_G(x_1)$ of $x_1$ and $N_G(x_2)$ of $x_2$ in $G$ except $x_1,x_2$ themselves, that is,
$N_{G'}(x_0)=(N_G(x_1)\cup N_G(x_2))\setminus\{x_1,x_2\}$, and
\item 
the red neighbours of $x_0$, denoted here by $N_{G'}^r(x_0)$, inherit all red neighbours of $x_1$ and of $x_2$ and add those in $N_G(x_1)\Delta N_G(x_2)$,
that is, $N_{G'}^r(x_0)=\big(N_{G}^r(x_1)\cup N_G^r(x_2)\cup(N_G(x_1)\Delta N_G(x_2))\big)\setminus\{x_1,x_2\}$, where $\Delta$ denotes the~symmetric set difference.
\end{itemize}
A \emph{contraction sequence} of a trigraph $G$ is a sequence of successive contractions turning~$G$  into a single vertex,
and its \emph{width} $d$ is the maximum red degree of any vertex in any trigraph of the sequence.
We also then say that it is a $d$-contraction sequence of~$G$.
The \emph{twin-width} of a trigraph $G$ is the minimum width over all possible contraction sequences of~$G$.
In other words, a graph has twin-width at most $d$ if and only if it admits a $d$-contraction sequence.

\smallskip
After the first implicit (and astronomical) upper bounds on the twin-width of planar graphs, e.g.~\cite{DBLP:journals/jacm/BonnetKTW22},
we have seen a stream of improving explicit bounds 
\cite{DBLP:journals/corr/abs-2202-11858,DBLP:conf/wg/JacobP22,DBLP:conf/isaac/BekosLH022,DBLP:journals/corr/abs-2205.05378},
culminating with the current best published upper bound of~$8$ by Hlin\v{e}n\'y and Jedelsk\'y~\cite{DBLP:conf/icalp/HlinenyJ23}.
This is complemented with a nearly matching lower bound of~$7$ by Kr\'al' and Lamaison in \cite{DBLP:journals/corr/abs-2209-11537}.
The right maximum value ($7$ or $8$?) is open, but the recent research of Jedelsk\'y~\cite{Jedelsky2024thesis}
strongly indicates that $7$ is the right answer, but this claim likely requires a computer-assisted proof.

It comes without surprise that the gradually improving upper bounds have required stronger and more involved arguments,
and the best ones are not easy to read for non-experts.
In this paper, we take the opposite route; we give a slightly worse bound with a self-contained proof which is as short
and simple as possible with the current knowledge:
\begin{theorem}\label{thm:twwplanar}
The twin-width of any simple planar graph is at most~$11$.
\end{theorem}

\section{Layered Skeletal Trigraphs}\label{sec:tools}

We use standard terminology of graph theory, and assume every graph to be simple (without loops and multiple edges).
We will mainly deal with planar graphs.
Recall that in $2$-connected planar graphs, every face is bounded by a cycle.
A {\em BFS tree} of a graph $G$ is a spanning tree defined by a run of the breadth-first-search algorithm on~$G$.

For a (tri)graph $G$, an ordered partition $\ca L=(L_0,L_1,\ldots)$ of $V(G)$ is called a {\em layering of~$G$} if,
for every edge $\{v,w\}$ of $G$ with $v\in L_i$ and $w\in L_j$, we have $|i-j|\leq 1$.
For example, every BFS tree $T\subseteq G$ with the root $r$ naturally defines a layering;
$L_0=\{r\}$, and $L_i$ for $i>0$ consisting of all vertices of $G$ at graph distance $i$ from~$r$.

If $T\subseteq G$ is a rooted tree (e.g., a BFS tree), 
a path $P\subseteq G$ is called {\em$T$-vertical} if $P\subseteq T$ is a subpath of some leaf-to-root path of~$T$.
A cycle $C\subseteq G$ is called {\em$T$-wrapped} if there exists an edge $e\in E(C)$, such that
$C-e$ is not $T$-vertical and $C-e$ is the union of two $T$-vertical paths intersecting in one vertex $u\in V(C)$.
Note that such $u$ must be unique -- it is the vertex of $C$ closest to the root of~$T$, and we call $u$ the {\em sink of~$C$}.
Moreover, observe that both edges incident to $u$ in $C$ belong to~$T$ as well.

Our goal, regarding \Cref{thm:twwplanar}, is to recursively construct $11$-contraction sequences of all planar graphs.
In order to formulate it, we now introduce our key concept -- of a ``splendid layered skeletal trigraph''.

\begin{definition}[Skeletal trigraph]\label{def:skeletal}\rm
Let $H$ be a trigraph and $S\subseteq H$ a $2$-connected nonempty planar subgraph such that 
all edges of $H$ induced by $V(S)$ are black (note; including the edges not in $E(S)$).
Fix a plane embedding of $S$, and call $S$ a {\em plane skeleton of~$H$}.
Further, call a mapping of connected components of $H-V(S)$ to faces of $S$ a {\em face assignment} of $H$ in~$S$,
if every connected component $H_0$ of $H-V(S)$ is assigned to a face $\phi$ of $S$ such that all neighbours of $H_0$ in $V(S)$ belong to~$\phi$.
Denote by $U_\phi$ the union of the vertex sets of all components assigned to~$\phi$ in this assignment.

If $H$ and $S$ conform to the previous conditions and there exists a face assignment of $H$ in~$S$, or if $S$ is the empty graph,
then we call $(H,S)$ a {\em skeletal trigraph}.
Moreover, if $\ca L$ is a layering of~$H$, then $(H,S,\ca L)$ is a {\em layered skeletal trigraph}.
\end{definition}

\begin{definition}[Splendid layered skeletal trigraph]\label{def:splendid}\rm
Consider a layered skeletal trigraph $(H,S,\ca L)$ with nonempty $S$ as in \Cref{def:skeletal}, and a face $\phi$ of~$S$.
We say that $\phi$ is {\em blank} if $U_\phi=\emptyset$ (i.e., if no connected component of $H-V(S)$ is assigned to~$\phi$),
and that $\phi$ is {\em$k$-reduced} if~$|U_\phi\cap L_i|\leq k$ holds for every layer $L_i\in\ca L$.
%
A layered skeletal trigraph $(H,S,\ca L)$ is {\em splendid} if 
either $S=\emptyset$ and $|V(H)\cap L_i|\leq4$ holds for all $L_i\in\ca L$ (i.e., whole $H$ is $4$-reduced),
or $S\not=\emptyset$ and all four following conditions are satisfied:
\begin{enumerate}[a)]\parskip 2pt
\item Every blank face of the plane skeleton $S$ is a triangle.
All non-blank faces of $S$, except at most one, are $1$-reduced, and the possible remaining face of $S$ is $3$-reduced.
\label{it:onerich}

\item There exists a BFS tree $T\subseteq S$ of the skeleton $S$ such that
the layering defined by $T$ in $S$ is equal to the restriction of $\cal L$ to~$V(S)$,
and that the facial cycle of every non-blank face $\phi$ of $S$ is $T$-wrapped.
\label{it:Tpropert}

\item For every non-blank face $\phi$ of $S$ with the facial cycle $C$, and $u$ being the sink of $C$, the following holds.
If $u\in L_i\in\ca L$, then all vertices of $U_\phi\cup V(C-u)$ belong to $L_{i+1}\cup L_{i+2}\cup\ldots$,
and there is a black edge in $H$ (but {\em no} red edge) from $u$ to each vertex of $U_\phi\cap L_{i+1}$.
\label{it:nosinkred}

\item Assume $\phi$ is a face of $S$ with the facial cycle $C$ which is not $1$-reduced, and $i$ is such that $L_i\in\ca L$.
Then every vertex $v$ in $X:=(U_\phi\cup V(C))\cap L_i$ has in~$H$ at most 
$3$ red edges into other vertices of $X$ and at most $4$ red edges into $U_\phi\cap(L_{i-1}\cup L_{i+1})$
(where, possibly,~$L_{-1}=\emptyset$).
Moreover, if $|U_\phi\cap L_{i+1}|>1$, then $v\in X$ has at most $2$ red edges into $U_\phi\cap L_{i-1}$.
\label{it:tworedin}
\end{enumerate}
\end{definition}

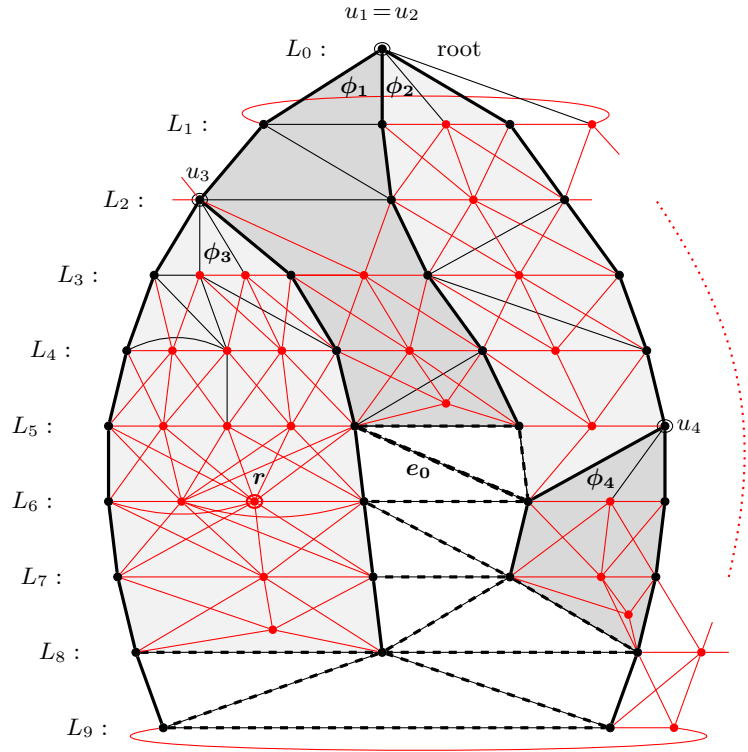
\begin{figure}[t]
\vspace*{-5ex}$$\quad
\begin{tikzpicture}[yscale=1.0, xscale=1.2]
\small
\tikzstyle{every node}=[draw, shape=circle, minimum size=3pt,inner sep=0pt, fill=black]
\node[label=right:\qquad root, label=left:$L_0:\qquad$, label=above:{$\!\!u_1\!=\!u_2\!\!$}] at (0,9) (l0) {};
\node[label=left:$L_1:\qquad$] at (-1.3,8) (l1) {}; \node at (0,8) (q1) {}; \node at (1.4,8) (r1) {};
\node[label=left:$L_2:\qquad$, label=above:{$\>u_3~$}] at (-2,7) (l2) {}; \node at (0.1,7) (q2) {}; \node at (2,7) (r2) {};
\node[label=left:$L_3:\qquad$] at (-2.5,6) (l3) {}; \node at (-1.0,6) (p3) {}; \node at (0.5,6) (q3) {}; \node at (2.6,6) (r3) {};
\node[label=left:$L_4:\qquad$] at (-2.8,5) (l4) {}; \node at (-0.5,5) (p4) {}; \node at (1.1,5) (q4) {}; \node at (2.9,5) (r4) {};
\node[label=left:$L_5:\qquad$] at (-3,4) (l5) {}; \node at (-0.3,4) (p5) {}; \node at (1.5,4) (q5) {}; \node[label=right:$\>u_4$] at (3.1,4) (r5) {};
\node[label=left:$L_6:\qquad$] at (-3,3) (l6) {}; \node at (-0.2,3) (p6) {}; \node at (1.6,3) (q6) {}; \node at (3.1,3) (r6) {};
\node[label=left:$L_7:\qquad$] at (-2.9,2) (l7) {}; \node at (-0.1,2) (p7) {}; \node at (1.4,2) (q7) {}; \node at (3,2) (r7) {};
\node[label=left:$L_8:\qquad$] at (-2.7,1) (l8) {}; \node at (0,1) (p8) {}; \node at (2.8,1) (r8) {};
\node[label=left:$L_9:\qquad$] at (-2.4,0) (l9) {}; \node at (2.5,0) (r9) {};
\draw (l0) circle (2.4pt); \draw (l2) circle (2.4pt); \draw (r5) circle (2.4pt);
\draw (q1)--(l1)--(q2)--(l2);
\draw (p5)--(q5)--(q6)--(p5)--(q4);
\draw (q6)--(p6)--(q7)--(p7) (p8)--(q7)--(r8)--(p8);
\draw (l9)--(r9)--(p8)--(l9) (l8)--(p8);
\draw (r4)--(q3)--(r2);
\draw[red] (2.3,8) to[out=20, in=160] (l1);
\begin{scope}[on background layer]
\draw[fill=white!95!black] (0,9)--(0,8)--(0.1,7)--(0.5,6)--(1.1,5)--(1.5,4)--(1.6,3)--(3.1,4)--(2.9,5)--(2.6,6)--(2,7)--(1.4,8);
\draw[fill=white!85!black] (0,9)--(0,8)--(0.1,7)--(0.5,6)--(1.1,5)--(1.5,4)--(-0.3,4)--(-0.5,5)--(-1.0,6)--(-2,7)--(-1.3,8);
\draw[fill=white!95!black] (0,1)--(-0.3,4)--(-0.5,5)--(-1.0,6)--(-2,7)--(-2.5,6)--(-2.8,5)--(-3,4)--(-3,3)--(-2.9,2)--(-2.7,1);
\draw[fill=white!85!black] (3.1,4)--(3.1,3)--(3,2)--(2.8,1)--(1.4,2)--(1.6,3);
\end{scope}
\tikzstyle{every path}=[very thick]
\draw (l0)--(l1)--(l2)--(l3)--(l4)--(l5)--(l6)--(l7)--(l8)--(l9);
\draw (l0)--(r1)--(r2)--(r3)--(r4)--(r5)--(r6)--(r7)--(r8)--(r9);
\draw (l2)--(p3)--(p4)--(p5)--(p6)--(p7)--(p8);
\draw (l0)--(q1)--(q2)--(q3)--(q4)--(q5);
\draw (r5)--(q6)--(q7);
\draw[dashed] (l9)--(r9)--(p8)--(l9) (l8)--(p8);
\draw[dashed] (p5)--(q5)--(q6)--(p5);
\draw[dashed] (q6)--(p6)--(q7)--(p7) (p8)--(q7)--(r8)--(p8);
\draw[ultra thick, dashed] (p5)--(q6);
\tikzstyle{every node}=[draw, shape=circle, minimum size=2pt,inner sep=0pt, color=red, fill=red]
\node at (2.3,8) (o1) {}; \node at (3.5,1) (o2) {}; \node at (3.2,0) (o3) {};
\node at (0.7,8) (z1) {}; \node at (1,7) (z2) {}; \node at (1.5,6) (z3) {};
\node at (1.9,5) (z4) {}; \node at (2.3,4) (z5) {};
\node at (-0.2,6) (y1) {}; \node at (0.3,5) (y2) {}; \node at (0.7,4.3) (y3) {};
\node at (2.5,3) (t1) {}; \node at (2.4,2) (t2) {}; \node at (2.7,1.5) (t3) {};
\node at (-2,6) (x1) {}; \node at (-1.5,6) (x2) {};
\node at (-2.3,5) (x3) {}; \node at (-1.7,5) (x4) {}; \node at (-1.1,5) (x5) {};
\node at (-2.4,4) (x6) {}; \node at (-1.7,4) (x7) {}; \node at (-1.0,4) (x8) {};
\node at (-2.2,3) (x9) {}; \node at (-1.4,3) (x10) {};
\node at (-1.3,2) (x11) {}; \node at (-1.2,1.3) (x12) {};
\tikzstyle{every path}=[]
\draw (l0)--(z1) (l0)--(o1) (r5)--(t1) (x7)--(x4)--(x1)--(l2)--(x2);
\draw (l3)--(x1) (l3)--(x4) (x1)--(p4) (l4) to[bend left] (x4);
\tikzstyle{every path}=[color=red]
\draw (r1)--(o1) (r2)--(o1)-- +(0.3,-0.4);
\draw (r8)--(o2)-- +(0.3,0) (r9)--(o2)-- +(0.12,0.4) (r7)--(o2)--(o3);
\draw (r9)--(o3)--(r8);
\draw (o3) to[out=350, in=190] (l9);
\draw[thick, dotted] (3.8,2) to[bend right=20] (3,7);
\draw (l2)-- +(-0.3,0) (l2)-- +(-0.2,0.3) (r2)-- +(0.3,0);
\draw (z1)--(z2)--(z3)--(z4)--(z5)--(q6);
\draw (q1)--(z1)--(r1)--(z2)--(q1);
\draw (q2)--(r2) (q2)--(z1)--(r2)--(z3)--(q2);
\draw (q3)--(r3) (q3)--(z2)--(r3)--(z4)--(q3);
\draw (q4)--(r4) (q4)--(z3)--(r4)--(z5)--(q4);
\draw (q5)--(r5);
\draw (l2)--(y1)--(y2)--(y3)--(p5);
\draw (p3)--(q3) (q3)--(y2)--(p3)--(y1)--(q2);
\draw (p4)--(q4) (p4)--(y1)--(q4)--(y3)--(p4);
\draw (y3)--(q5)--(y2)--(p5);
\draw (t1)--(t2)--(t3)--(r8) (t2)--(r8);
\draw (q6)--(r6) (q6)--(t2)--(r6);
\draw (q7)--(r7) (q7)--(t1)--(r7)--(t3)--(q7);
\draw (x1)--(p3) (l4)--(p4) (l5)--(p5) (l6)--(p6) (l7)--(p7);
\draw (x2)--(x5)--(x8)--(x10)--(x11)--(x12);
\draw (x1)--(x3)--(x6)--(x9)--(x11);
\draw (l3)--(x3)--(l5) (x1)--(x5)--(p3);
\draw (x4)--(x2)--(p4) (l4)--(x6)--(x4)--(x8)--(p4);
\draw (x3)--(x7)--(x5)--(p5);
\draw (l5)--(x9)--(x7)--(x10)--(p5) (l5)--(x10)--(x6) (x8)--(x9)--(p5);
\draw (x6)--(l6)--(x11) (x7)--(p6)--(x11) (x8)--(p6);
\draw (l6) to[bend right=20] (x10) (x9) to[bend right=20] (p6);
\draw (l7)--(x9)--(p7)--(x10)--(l7);
\draw (l8)--(x11)--(p8)--(x12)--(l8) (l7)--(x12)--(p7);
\draw[thick] (x10) circle (2.4pt);
\tikzstyle{every node}=[draw=none, color=black]\boldmath
\node at (-1.8,6.3) {$\phi_3$};
\node at (-0.3,8.5) {$\phi_1$};
\node at (0.2,8.5) {$\phi_2$};
\node at (2.4,3.3) {$\phi_4$};
\node at (-1.35,3.3) {$r$};
\node at (0.4,3.4) {$e_0$};
\end{tikzpicture}
$$
\caption{A picture of a splendid layered skeletal trigraph $(H,S,\ca L)$, in which the skeleton $S$ is depicted with black vertices
	and thick black edges such that the associated BFS tree $T\subseteq S$ is drawn with thick solid edges and the edges of $E(S)\setminus E(T)$ are thick dashed.
	$T$ has its root at the top and its (ten) BFS layers are organized horizontally in the picture.
	There are four bounded non-blank faces in $S$, denoted by $\phi_1,\phi_2,\phi_3,\phi_4$ (with corresponding sinks $u_1,u_2,u_3,u_4$), and emphasized with gray shade.
	The unbounded face of $S$ is also non-blank, but it is only sketched in the picture.
	There is one non-$1$-reduced face in $(H,S)$, namely $\phi_3$, and it contains a red vertex $r$ (emphasized with a circle around) that achieves the
	maximum red degree $11$ allowed by \Cref{def:splendid}.}
\label{fig:splendid}
\end{figure}

\Cref{def:splendid} is illustrated, with comments, in \Cref{fig:splendid}.
%
In regard of the definition we stress that the subgraph $H_\phi\subseteq H$ of an $S$-face $\phi$ induced by $U_\phi\cup V(C)$ need not be planar 
(since non-planarity may easily be introduced by contractions), and some vertices of $H_\phi$ may actually belong to layers
of $\ca L$ which are higher than the maximum layer intersecting $V(C)$.

The following simple proof also approachably illustrates \Cref{def:splendid} (especially~\ref{def:splendid}.\ref{it:tworedin}).
\begin{lemma}\label{lem:atmost11}
Every splendid layered skeletal trigraph has maximum red degree at most~$11$.
\end{lemma}
\begin{proof}
Let $(H,S,\ca L)$ be a splendid layered skeletal trigraph.
By \Cref{def:skeletal}, every red edge of $H$ must have one or both ends in~$V(H)\setminus V(S)$,
and must be within the same or consecutive layers of~$\ca L$.
Hence for $S=\emptyset$ we immediately get an upper bound of $4+3+4=11$ on the red degree.
For the rest we assume $S\not=\emptyset$.

Consider a vertex $v\in V(S)$. Then $v$ has no red edges into vertices of $S$ by \Cref{def:skeletal}.
If $\phi$ is a face of $S$ incident to~$v$, then,
by \Cref{def:splendid}.\ref{it:nosinkred}, $v$ can have red edges into $U_\phi$ only if $v$ is not the sink of~$\phi$.
Let $r$ be the root of the BFS tree $T\subseteq S$ by \Cref{def:splendid}.\ref{it:Tpropert}.
Then $v=r$ is the sink of every incident face, and for every $v\not=r$ we have only two faces of $S$ such that
$v$ is not the sink of them -- these are the two faces adjacent to the parental edge of $v$ in $T$ (the edge pointing towards the root~$r$).
Let these two faces of $v$ be $\phi_1$ and $\phi_2$.
By \Cref{def:splendid}.\ref{it:onerich}, we may up to symmetry assume that $\phi_1$ is $1$-reduced, and then $v$ can have at most
$1+1+1=3$ red neighbours in $U_{\phi_1}$. The same is true if $\phi_2$ is $1$-reduced.
Moreover, by \Cref{def:splendid}.\ref{it:tworedin}, $v$~can have at most $3+4=7$ red neighbours in $U_{\phi_2}$ even if $\phi_2$ is not $1$-reduced.
Altogether,~$3+3<3+7<11$.

Consider now $v\in V(H)\setminus V(S)$, and the (non-blank) face $\phi$ such that $v\in U_\phi$.
Let $C\subseteq S$ be the cycle bounding~$\phi$; then all red neighbours of $v$ belong to $U_\phi\cup V(C)$ by \Cref{def:skeletal}.
Observe also that $|V(C)\cap L_i|\leq2$ for all $L_i\in\ca L$ by \Cref{def:splendid}.\ref{it:Tpropert}.
If $\phi$ is $1$-reduced, then the claim easily holds.
Otherwise, there are at most $3+4=7$ red neighbours of $v$ accounted for by \Cref{def:splendid}.\ref{it:tworedin}, 
and additional at most $4$ in the set $V(C)\cap(L_{i-1}\cup L_{i+1})$, where $v\in L_i\in\ca L$. 
Again at most~$7+4=11$ altogether.
\end{proof}

The main statement of the paper follows.

\begin{lemma}\label{lem:induction}
Every splendid layered skeletal trigraph admits an $11$-contraction sequence.
\end{lemma}
While we leave the inductive proof of this lemma to the next section, we show how it implies our main result.

\begin{proof}[Proof of \Cref{thm:twwplanar}]
Given a planar graph $G$, we fix any plane embedding of~$G$.
We construct a plane triangulation $G^+$ from $G$ by adding new vertices to every non-triangular face of $G$ and connecting them inside each face and to vertices of this face.%
\footnote{If a face $\varphi$ is bounded by a cycle, we simply add one new vertex adjacent to all boundary vertices, but that could violate
simplicity of $G^+$ if the boundary of $\varphi$ is not a cycle. One can check that the proof here works with non-simple graphs as well (only the
definition of twin-width needs a simple graph), but it is cleaner to add into $\varphi$ one more new vertex for every repetition of a vertex 
on the boundary of~$\varphi$ and keep $G^+$ simple.}
Then $G^+$ is $2$-connected.
Choosing an arbitrary BFS tree of $G^+$, we take the layering $\ca L=(L_0,L_1,\ldots)$ of $G^+$ naturally defined by~$T$.
Then, trivially, $(G^+,G^+,\ca L)$ is a splendid layered skeletal trigraph, and hence $G^+$ admits an $11$-contraction sequence by \Cref{lem:induction}.
Restricting this sequence only to the contractions of pairs from $V(G)$ we, again trivially,
obtain an $11$-contraction sequence of~$G$.
\end{proof}

\section{Proof of \Cref{lem:induction}}\label{sec:induction}

Our proof starts with an auxiliary claim.
\begin{lemma}\label{lem:bound1leaf}
Let $G$ be a $2$-connected plane graph, and $T\subseteq G$ a BFS tree of~$G$.
Assume~$T$ that has at least $3$ leaves, and that for every facial cycle $C$ of $G$, we have $|E(C)\setminus E(T)|=1$ or $C$ is a triangle.
Then there exists an edge $e\in E(G)\setminus E(T)$ such that, for the unique cycle $D_e\subseteq T+e$, one of the two faces of $D_e$
contains (in its strict interior) precisely one leaf of $T$ and not the root of $T$.
\end{lemma}

\begin{proof}
Picture $G$ in the plane such that the root of $T$ is on the unbounded face.
Then no bounded face of $G$ contains the root in its interior.
There exists an edge $e\in E(G)\setminus E(T)$ such that the interior of the bounded face of the cycle $D_e$ 
in the considered plane drawing of~$G$ contains some leaf of~$T$; one can choose the outer face boundary of $G$ by the assumptions of the lemma.
Among all such possible choices of $e$, we select one such that the interior of $D_e$ contains the least number of vertices of $G$.
%
%
By means of contradiction, we assume that this $e$ is not a sought solution, meaning that the interior of $D_e$      
contains at least two leaves of~$T$.

We have $e\not\in E(T)$ and consider the face $\sigma$ of $G$ in the interior of $D_e$ and adjacent to~$e$.
Since the interior of $D_e$ contains a vertex and $G$ is $2$-connected,
it cannot happen that $\sigma$ is bounded by $D_e$, and so $\sigma$ has another edge not from $T$.
Hence, by the assumption on $G$, the face $\sigma$ is bounded by a triangle $C_0$ such that $E(C_0)=\{e,e_1,e_2\}$.
Moreover, some of the at least two leaves of $T$ in the interior of $D_e$ is not on $C_0$.
If, up to symmetry, $e_1\in E(T)$, then we can choose $e_2$ instead of $e$ which contradicts minimality of our choice of~$e$.
Otherwise, for at least one of the choices of $e_i$, $i\in\{1,2\}$, instead of $e$, the cycle $D_{e_i}\subseteq T+e_i$
encloses at least one leaf of $T$ in the interior, and we again have a contradiction to our minimal choice of~$e$ above.
\end{proof}

We proceed to the {\sf\bfseries\textcolor{lipicsGray}{Proof of \Cref{lem:induction}}}, considering a splendid layered skeletal trigraph $(H,S,\ca L)$.
For start, the maximum red degree of $H$ is at most $11$ by \Cref{lem:atmost11}.
For the rest of a sought $11$-contraction sequence of $H$, we proceed by induction on $|V(H)|+|V(S)|$.

If the skeleton is $S=\emptyset$, then we pick the largest index $i$ such that $V(H)\cap L_i\not=\emptyset$.
If $|V(H)\cap L_i|>1$, we contract any two vertices in $V(H)\cap L_i$,
and if $V(H)\cap L_i=\{x\}$ and $i\geq1$, we contract $x$ with any vertex of $V(H)\cap L_{i-1}$.
In both cases, the contraction results in a splendid layered skeletal trigraph, again with $S=\emptyset$,
and so we may finish by induction.
Otherwise, that is for $V(H)\cap L_i=\{x\}$ and $i=0$, we are done as $V(H)=\{x\}$.

Thus, we may assume that $S\neq\emptyset$.
If all faces of $S$ are $1$-reduced or blank, and the BFS tree $T\subseteq S$ from \Cref{def:splendid}.\ref{it:Tpropert} 
has at most $2$ leaves, we get that $T$ consists of at most two $T$-vertical paths, and that $S$ has at most two non-blank faces by \Cref{def:splendid}.\ref{it:Tpropert}.
Since the two faces are $1$-reduced, every layer of $\ca L$ contains at most $1+1+2=4$ vertices.
So, $(H,S'=\emptyset,\ca L)$ is also a splendid layered skeletal trigraph and we continue as in the previous paragraph.

Thus, we may assume that the BFS tree $T\subseteq S$ has at least $3$ leaves 
or the skeleton $S\not=\emptyset$ has a face which is not $1$-reduced.
We have two cases.

\vspace*{-2ex}
\subparagraph{Case 1.} The skeleton $S$ has all faces $1$-reduced (some or all may be blank).

Then $T$ has at least $3$ leaves.
We apply \Cref{lem:bound1leaf} and get the edge $e$ and cycle $D_e\subseteq T+e\subseteq H$
such that in the interior of $D_e$ there is precisely one leaf $x$ of $T$.
Let $Q$ be the maximal $T$-vertical path starting in $x$ and not hitting~$D_e$.
Importantly, all vertices of $S$ in the interior of $D_e$ must lie on $Q$, or there were another leaf of $T$ there since $T$ is spanning.

The interior of~$D_e$ contains at most two non-blank faces of $S$; this is since
there are at most two available sinks for the faces -- the sink $u$ of $D_e$ and the vertex $w$ of $D_e$
to which $Q$ is adjacent in $T$ (these two may coincide $u=w$, and still be the sinks of two faces).
Let these faces of $S$ be $\phi_1$ and $\phi_2$, and note that either $\phi_2$ or both $\phi_1,\phi_2$ may possibly be blank.
The considered case can be illustrated in \Cref{fig:splendid} (ignoring for now that the face $\phi_3$ is not $1$-reduced) with the
edge $e=e_0$ chosen by \Cref{lem:bound1leaf}, such that the cycle $D_{e_0}$ with the sink $u_1=u_2$
encloses two $1$-reduced faces $\phi_1,\phi_2$ and one blank triangular face.
In general, there can be more than one blank faces of the skeleton $S$ enclosed by $D_{e}$.

We set $S':=S-V(Q)$ and consider the layered skeletal trigraph $(H,S',\ca L)$ with the (new) non-blank face $\phi$ bounded by~$D_e$.
The first step is to observe that $\phi$ is a $3$-reduced face of $(H,S',\ca L)$, as required by \Cref{def:splendid}.\ref{it:onerich}:
since $U_\phi=U_{\phi_1}\cup U_{\phi_2}\cup V(Q)$, for every~$L_j\in\ca L$, the set
$U_\phi\cap L_j$ has at most one vertex from $Q$, and at most one from each of $U_{\phi_k}\cap L_j$, $k=1,2$ since $\phi_k$ was $1$-reduced
(or possibly blank which means $U_{\phi_k}=\emptyset$).
The conditions of \Cref{def:splendid}.\ref{it:Tpropert} are fulfilled (with the restriction of~$T$) by the choice of $e$ in \Cref{lem:bound1leaf},
and the first part of \Cref{def:splendid}.\ref{it:nosinkred} follows analogously since the root of $T$ is on or outside of~$D_e$.
For the second part of \Cref{def:splendid}.\ref{it:nosinkred}, take the sink $u\in L_i$ of $\phi$ and a vertex $x\in U_\phi\cap L_{i+1}$.
If $x\in U_{\phi_1}\cup U_{\phi_2}$, then there is a black edge $ux$ in $H$ 
by \Cref{def:splendid}.\ref{it:nosinkred} applied to~$(H,S,\ca L)$, and if $x\in V(Q)$, then $ux$ is a black edge in~$T$.

It remains to check the conditions of \Cref{def:splendid}.\ref{it:tworedin} for $(H,S',\ca L)$.
Pick a vertex $x\in(U_\phi\cup V(D_e))\cap L_i$.
If $x\in V(Q)\cup V(D_e)$, then $x$ can have red neighbours in $H$ only in $U_{\phi_1}\cup U_{\phi_2}$ by the skeletal trigraph~$(H,S)$,
which means at most $2$ red edges into every layer of~$\ca L$ there (and so also into each of the sets~$X$ and $U_\phi\cap L_{i-1}$), as required.
If $x\in U_\phi\setminus V(Q)$ and, up to symmetry, $x\in U_{\phi_1}$, then red edges from $x$ may only end in 
$Y:= U_{\phi_1}\cup V(Q)\cup V(D_e)$, again by the skeletal trigraph~$(H,S)$.
By cardinalities of the respective sets,
there are at most $3$ red neighbours in $(Y\setminus\{x\})\cap L_i\subseteq V(Q)\cup V(D_e)$, and at most $2$ red neighbours 
in either of $(U_{\phi_1}\cup V(Q))\cap L_{i-1}$ or $(U_{\phi_1}\cup V(Q))\cap L_{i+1}$.
In particular, $v$ has (unconditionally) at most $2$ red edges into $U_\phi\cap L_{i-1}$ since those cannot end in $U_{\phi_2}$ by the skeletal trigraph~$(H,S)$.

Therefore, $(H,S',\ca L)$ is splendid, and we finish by induction with it.

\vspace*{-2ex}
\subparagraph{Case 2.} The skeleton $S$ has a face $\phi$ which is not $1$-reduced.

By \Cref{def:splendid}.\ref{it:onerich},
$\phi$ is $3$-reduced, and let $j$ be the largest index such that $|U_\phi\cap L_j|>1$ for~$L_j\in\ca L$.
We contract any two vertices $v,w\in U_\phi\cap L_j$ in $H$ into a new vertex $t$, creating a layered skeletal trigraph
$(H',S,\ca L')$ where $H'$ results from $H$ by this contraction and the layering $\ca L'$ is naturally inherited from~$\ca L$
(i.e., $t\in L'_j\in\ca L'$). 
Denoting by $U'_\phi=(U_\phi\setminus\{v,w\})\cup\{t\}$ the set that stems from $U_\phi$, we get $|U'_\phi\cap L'_{j}|\leq2$.
For an illustration, see the face $\phi=\phi_3$ in \Cref{fig:splendid} in which the trigraph $H'$ (as called $H$ there) resulted by a contraction
of two vertices from $U_{\phi_3}\cap L_6$ into the emphasized vertex~$r$.

We are going to prove that $(H',S,\ca L')$ satisfies the conditions of \Cref{def:splendid}, and then apply induction.
Validity of \Cref{def:splendid}.\,\ref{it:onerich},\ref{it:Tpropert} follows directly from the same conditions for $(H,S,\ca L)$.
As for \Cref{def:splendid}.\ref{it:nosinkred}, the situation changes only if the sink of $\phi$ is $u\in L_i$ where $i=j-1$.
Then \ref{def:splendid}.\ref{it:nosinkred} remains valid in $H'$ since each of $v$ and $w$ had a black edge to $u$ in $H$.

Checking preservation of the conditions of \Cref{def:splendid}.\ref{it:tworedin} is relevant whenever $j-1\leq i\leq j+1$.
For $i=j+1$ we already have $|U'_\phi\cap L'_{i+1}|\leq1$ and $|U'_\phi\cap L'_{i}|\leq1$ by the definition of~$j$,
and $|U'_\phi\cap L'_{i-1}|\leq2$ as noted above.
So, the conditions are satisfied in~$H'$ as well.
For $i=j$ we again have $|U'_\phi\cap L'_{i+1}|\leq1$, and so in~$H'$, any vertex $v$ from $X=(U'_\phi\cup V(C))\cap L'_i$ can have 
only $1$ neighbour in $U'_\phi\cap L'_{i+1}$ and at most $|U'_\phi\cap L'_{i-1}|\leq3$ neighbours in $U'_\phi\cap L'_{i-1}$, as claimed.
The number of neighbours of $v$ in the rest of $X$ is bounded by $2+|U'_\phi\cap L'_{i}|-1\leq2+2-1=3$.
For $i=j-1$, from the assumption $|U_\phi\cap L_{j}|>1$ and \Cref{def:splendid}.\ref{it:tworedin},
we get that there have been at most $2$ red neighbours of $v\in X$ in $U_\phi\cap L_{i-1}$ 
and these are not affected by the contraction in $L_{i+1}$. Hence the same is true in~$H'$ and~$U'_\phi$.
Furthermore, there are at most $|U'_\phi\cap L'_{i+1}|\leq2$ neighbours of $v$ in $U'_\phi\cap L'_{i+1}$, again as claimed.

Therefore, $(H',S,\ca L')$ is splendid, and we again finish by induction with it.
The whole proof of \Cref{lem:induction} is done.
\hfill$\Box$

\section{Conclusion}

We have provided a short self-contained proof of \Cref{thm:twwplanar}.
While the proved bound is not the best currently possible, 
the proof given here is way much simpler than those in \cite{DBLP:journals/corr/abs-2205.05378,DBLP:conf/icalp/HlinenyJ23}.
While sacrificing a bit of simplicity of the given proof, we can also give a better upper bound
of~$9$ (thus matching \cite{DBLP:journals/corr/abs-2205.05378}), but
we are so far not sure whether a similarly simplified proof can be given for the upper bound of $8$ as in~\cite{DBLP:conf/icalp/HlinenyJ23}.

The proof of \Cref{thm:twwplanar} is constructive, and it is not difficult to derive a polynomial-time algorithm
for a construction of the claimed $11$-contraction sequence from it; in particular, the edge $e$ of \Cref{lem:bound1leaf}
can be found by testing all relevant edges.
However, the (more complex) proof in~\cite{DBLP:conf/icalp/HlinenyJ23} comes along with a significantly simpler
linear-time algorithm for a construction of an $8$-contraction sequence, and so we skip the algorithmic aspects in this paper.

\bibliography{tww}

\end{document}